\input amstex
\magnification=\magstep0
\documentstyle{amsppt}
\pagewidth{6.2in}
\topmatter
\title
$q$-Ap\'ery Irrationality Proofs by $q$-WZ pairs
\endtitle
\author
Tewodros Amdeberhan and Doron Zeilberger
\endauthor
\affil Department of Mathematics, Temple University,
Philadelphia PA 19122, USA \\
tewodros\@math.temple.edu, zeilberg\@math.temple.edu
\endaffil
\abstract
Using WZ pairs, Ap\'ery-style proofs of the irrationality of the q-analogues of 
the Harmonic series and $Ln(2)$ are given. For the q-analogue of $Ln(2)$, 
this method produces an improved irrationality measure.
\endabstract
\endtopmatter
\def\({\left(}
\def\){\right)}

\document

\bf 0. Introduction: \rm \smallskip

Let us define the following q-analogues of the Harmonic series 
$\sum_{n=1}^{\infty}\frac 1n$ and Ln(2), respectively by:
$$\alignat2
h_q(1):&=\sum_{k=1}^{\infty}\frac 1{q^k-1} &&\qquad\text{(for $|q| > 1$),} 
\tag0.1
\endalignat$$
$$\alignat2
Ln_q(2):&=\sum_{n=1}^{\infty}\frac {(-1)^n}{q^n-1} &&\qquad\text{(for $|q| \neq 
0, 1$)}.\tag0.2
\endalignat$$

In 1948, Paul Erd\"os [E1] proved the irrationality of $h_2(1)$. Recently, 
Peter Borwein used Pad\'e approximation techniques [B1] and some complex analysis 
methods [B2] to prove the irrationality of both $h_q(1)$ and $Ln_q(2)$. Here we 
present a proof in the spirit of Ap\'ery's magnificent proof of the 
irrationality of $\zeta(3)$ [A], which was later delightfully accounted by Alf 
van der Poorten [P]. This method of proof gives favorable irrationality measure 
(=4.80) for $Ln_q(2)$ campared to the irrationality measure (=54.0) implied 
in [B1], [B2]. Further discussion of irrationality results for certain series 
is to be found in Erd\"os [E2]. 
\smallskip

We will assume familiarity with ref. [Z]. In particular,  \par
$\binom{n}k_q := \frac {(q)_n}{(q)_k (q)_{n-k}}$, where $(q)_0:=1$ and 
$(q)_n: = (1-q) \cdots (1-q^n),$ for $n \geq 1.$ \smallskip
N and K are \it forward shift operators \rm on $n$ and $k$, respectively.
\smallskip
$\Delta_n := N - 1$, $\Delta_k := K - 1$. \smallskip
A pair $(F(n,k),G(n,k))$ of discrete functions is called a \it q-WZ pair \rm
if:\par
1. $NF/F$, $KF/F$, $NG/G$ and $KG/G$ are all rational functions of $q^n$ and
$q^k$, and \par
2. $\Delta_n F = \Delta_k G.$\smallskip

Given such a pair $(F,G)$, then
$\omega = F(n,k) \delta k + G(n,k)\delta n$ is called a \it q-WZ 1-form. \rm

\pagebreak
\midspace{.3in}

\bf 1. A scheme for proving the irrationality of the \it q-harmonic series \bf 
$h_q(1)$:\rm \bigskip
The claims made in subsections \bf 1.1-1.5 \rm below were found using the Maple 
Package 
{\tt qEKHAD} accompanying [PWZ]. The relevant script substantiating our claims 
can be found in this paper's Web Pages.\smallskip

\bf 1.1. \rm The q-WZ 1-form $\omega$ is:
$$\omega = \frac {-1}{\binom {n+k+1}{k}_{q}(q)_{n+1}}
\{\delta k + \frac {q^{n+1}}{(q^{n+1}-1)}\delta n\}.$$

\bf 1.2. \rm The choice of the potential $c(n,k)$ is:
$$c(n,k)=\sum_{m=1}^{n}\frac {q^m}{(1-q^m)(q)_{m}} + \sum_{m=1}^{k}\frac 
1{(q^m-1)}\frac 1{\binom {n+m}m_{q}(q)_{n}}.$$

\smallskip

\bf 1.3. \rm The choice of the mollifier $b(n,k)$ is:
$$b(n,k)=(-1)^k q^{k(k+1)/2}\binom {n+k}k_{q}\binom nk_{q}.$$

\bf 1.4. \rm We define two sequences:
$$
a(n)=\sum_{k=0}^{n}c(n,k)b(n,k), \qquad\text{and} \qquad 
b(n)=\sum_{k=0}^{n}b(n,k).$$

\bf 1.5. \rm Introduce $L=y_2(n)N^2+y_1(n)N+y_0(n)$ and $B(n,k) = P_q^1(n,k)
b(n+1,k)$, 
where
$$
 A(n,k)=c(n,k)B(n,k)+\frac{(-1)^kq^{2n+3}}{q^{n+1}-1}\binom {n+1}k_q\frac 
{q^{\binom k2}}{(q)_{n+2}}P_q^2(n,k) \qquad\text{and}$$

$P_q^1(n,k)=-q\alpha_n^2\beta_k^{-1}(q^2\alpha_n+2q)+q\alpha_n^2(q^2\alpha_n^3+
2q(q+1)\alpha_n^2+3q\alpha_n-(q+1)-(\alpha_n+2)\beta_k)$

$$\multline
P_q^2(n,k)=q^2\alpha_n^2+q\alpha_n-2+\beta_k(q^2\alpha_n^5+q(2q+1)\alpha_n^4-2\
alpha_n^3-\alpha_n^3\beta_k-(2-q^{-1})\alpha_n^2\beta_k\\
-(3q+5)\alpha_n^2+2q^{-1}\alpha_n\beta_k+(q-1+2q^{-1})\alpha_n+(1+3q^{-1})),
\endmultline$$

$y_0(n)= q(\alpha_n-1)(q\alpha_n+2)$, ${ }$ $y_2(n)=(q\alpha_n-1)(\alpha_n+2)$, 
$ {}$ $\alpha_n=q^{n+1}$, $ { }$ $\beta_k=q^{k+1}$ and

\bigskip
$y_1(n)=q^3\alpha_n^5+2q^2(q+1)\alpha_n^4
+q^2\alpha_n^3-4q(q+1)\alpha_n^2
+(q^2-4q+1)\alpha_n+2(q+1).$

\bigskip
Then $$\alignat3
L(b(n,k))&=B(n,k)-B(n,k-1) &\qquad\text{and} &\qquad 
L(b(n,k)c(n,k))&=A(n,k)-A(n,k-1).\tag{$*$}\endalignat$$
Now, summing over $k$ in $(*)$ shows that both sequences $a(n)$ and $b(n)$ are 
solutions of $Lu(n) = 0$.

\pagebreak
\midspace{.4in}

\bf 1.6. \rm Set $b_n=b(n)$ and $a_n=a(n)$. Now, since $b_{n+1} > b_n$ and  $Lb_n = 0$, 
that is, \par $y_2(n)b_{n+2}+y_1(n)b_{n+1}+y_0(n)b_n = 0$, then asymptotically
 we have 
that
$$\frac {b_{n+2}}{b_{n+1}} = O\(\frac{y_1(n)}{y_2(n)}\) = O\(q^{3n+3}\).$$
Hence, 
$$b_n = O\(q^{\frac{3n^2}2}\).\tag1.6.1$$

On the other hand, $La_n = 0$ and $Lb_n = 0$ lead to the system of recurrence 
relations,
$$\alignat2
y_2(n)a_{n+2}+y_1(n)a_{n+1}+y_0(n)a_n = 0,&\qquad
y_2(n)b_{n+2}+y_1(n)b_{n+1}+y_0(n)b_n = 0.\tag1.6.2\endalignat$$

Multiplying out the first and the second equations in (1.6.2), respectively by 
$b_{n+2}$ and $a_{n+2}$, and subtracting we obtain
$$y_1(n)(a_{n+1}b_{n+2}-b_{n+1}a_{n+2}) = y_0(n)(a_{n+2}b_{n}-b_{n+2}a_{n}).$$

Rewriting this in the form 
$$\frac {a_{n+1}}{b_{n+1}}- \frac {a_{n+2}}{b_{n+2}} =
\frac{y_0(n)}{y_1(n)}\frac {b_n}{b_{n+1}}\(\frac{a_{n+2}}{b_{n+2}}- \frac 
{a_{n}}{b_{n}}\)$$

leads to the estimate

$$\left|\frac {a_{n+1}}{b_{n+1}}- \frac {a_{n+2}}{b_{n+2}}\right| \leq
\left|\frac{y_0(n)}{y_1(n)}\frac {b_n}{b_{n+1}}\(\frac{a_{n+2}}{b_{n+2}}- \frac 
{a_{n+1}}{b_{n+1}}\)\right|+
\left|\frac{y_0(n)}{y_1(n)}\frac {b_n}{b_{n+1}}\(\frac{a_{n+1}}{b_{n+1}}- \frac 
{a_{n}}{b_{n}}\)\right|,$$

which in turn yields 
$$\frac {a_{n+1}}{b_{n+1}}- \frac {a_{n}}{b_{n}} = O\(b_n^{-2}\).\tag1.6.3$$

Therefore,
$$h_q(1)-\frac {a_n}{b_n}=O\(b_n^{-2}\).\tag1.6.4$$
In particular, the sequence of rational numbers $\frac {a_n}{b_n}$ converges 
moderately quickly to $h_q(1)$.
\pagebreak
\midspace{.3in}

\bf 1.7. \rm For a given prime $p$, let $ord_p$k denote the exponent of $p$ 
in the prime expansion of $k$. Then we observe that
$$ord_p \binom nm_q \leq ord_p(q)_n - ord_p(q)_m.\tag1.7.1$$
\bf Note: \rm $$\binom{n+k}k_q\binom km_q = 
\binom{n+m}m_q\binom{n+k}{k-m}_q.\tag1.7.2$$

\bf Lemma 1: \rm The sequences 
$$
u_n=a_n(q)_{n+1}\prod_{s=[n/2]}^n(1-q^s) \qquad\text{and} \qquad z_n= 
b_n(q)_{n+1}\prod_{s=[n/2]}^n(1-q^s)$$
are polynomials in $q$ with integer coefficients, and moreover 
$$z_n = O\(q^{19n^2/8}\).\tag1.7.3$$
  
\bf Proof: \rm Applying (1.7.1) and (1.7.2), we can estimate the denominator of 
$u_n$ as: 
$$\align
ord_p\(\frac {(q^m-1)(q)_{n}\binom{n+m}m_q}{\binom{n+k}k_q}\)
&\leq ord_p\(\frac {(q^m-1)(q)_{n}\binom km_q}{\binom{n+k}{k-m}_q}\)\\
&\leq ord_p(q)_{n} + ord_p(q^m-1) + ord_p(q)_k - ord_p(q)_m\\
&\leq ord_p(q)_{n} + ord_p\prod_{s=[n/2]}^n(1-q^s) +  ord_p(q)_k - ord_p(q)_m\\
&\leq ord_p\((q)_{n}\prod_{s=[n/2]}^n(1-q^s)\),\endalign$$
since $m \leq k \leq n$. This proves the claim on $u_n$. And (1.7.3) follows 
from (1.6.1). The rest is trivial. \bigskip
 
\bf Lemma 2: \rm $h_q(1)$ - $\frac {u_n}{z_n}$ = $O\(\frac 
1{z_n^{1+\delta}}\)$; where $\delta =0.26316\dots$$>$ 0.\bigskip

\bf proof: \rm From (1.6.1), (1.6.4) and (1.7.3), we gather that
$$
h_q(1) - \frac{u_n}{z_n}=O\(b_n^{-2}\)
=O\(q^{-3n^2}\)
=O\(z_n^{-1-(5/19)}\).$$
\smallskip
Thus, we have proved:\smallskip
\bf Theorem 1: \rm If $|q|>1$ is an integer, $h_q(1)$ is irrational with 
irrationality measure 4.80. \smallskip

\bf Remark 1: \rm By invoking Theorem 7 ([Z], p.596) with $\omega$ as in 1.1, 
we obtain the series acceleration:
$$
h_q(1)=\sum_{n=1}^{\infty}\frac {q^n}{(1-q^n)(q)_n} \qquad\text{and}\qquad
h_q(1)=\sum_{n=1}^{\infty}\frac {1-q^n-q^{2n}}{(q^n-1)\binom {2n}n_q (q)_n}.$$

\pagebreak
\midspace{.3in}

\bf 2. A scheme for proving the irrationality of $Ln_q(2)$: \rm \bigskip

The claims made in subsections \bf 2.1-2.5 \rm below were found using the Maple 
Package 
{\tt qEKHAD} accompanying [PWZ]. The relevant script substantiating our claims 
can be found in this paper's Web Pages.\smallskip

\bf 2.1. \rm The qWZ 1-form $\omega$ is:
$$\omega = \frac {(-1)^k}{(1-q^{k+1})}\frac {(q)_n}{\binom 
{n+k+1}{k+1}_{q}(q^2)_{n}}
\{\delta k + \frac {q^{n+1}}{(1+q^{n+1})}\delta n\}.$$

\bf 2.2. \rm The choice of the potential $c(n,k)$ is:
$$c(n,k)=\sum_{m=1}^{n}\frac {q^m(q)_m}{(1-q^m)(q^2)_{m}} + \sum_{m=1}^{k}\frac 
{(-1)^{m-1}}{(1-q^m)}\frac {(q)_n}{\binom {n+m}m_{q}(q^2)_{n}}.$$

\smallskip

\bf 2.3. \rm The choice of the mollifier $b(n,k)$ is:
$$b(n,k)=q^{k(k+1)/2}\binom {n+k}k_{q}\binom nk_{q}.$$

\bf 2.4. \rm We define two sequences:
$$
a(n)=\sum_{k=0}^{n}c(n,k)b(n,k), \qquad\text{and} \qquad 
b(n)=\sum_{k=0}^{n}b(n,k).$$

\bf 2.5. \rm Introduce $L=y_2(n)N^2+y_1(n)N+y_0(n)$ and $B(n,k) = P_q^1(n,k)
b(n+1,k)$, 
where
$$
 A(n,k)=c(n,k)B(n,k)+\frac{(-1)^kq^{2n+3}}{1-q^{n+1}}\binom {n+1}k_q\frac 
{q^{\binom k2}(q)_{n+1}}{(q^2)_{n+1}}P_q^2(n,k)\qquad\text{and}$$

$$\multline
P_q^1(n,k)=q\alpha_n^2\bigl[q^3\alpha_n^5+q^2(1+q)\alpha_n^4+2q(1+q^2)\alpha_n^
3-(1-q+q^2)\alpha_n-3(1+q)\big]\\
+q\alpha_n^2\bigl[q\beta_k^{-1}(q^2\alpha_n^3+q(1+q)\alpha_n^2+(2-q)\alpha_n-2)
+(q\alpha_n^3+(q-1)\alpha_n^2+(2q-1)\alpha_n)\alpha_k-2\bigr] \endmultline$$
$$\multline
P_q^2(n,k)=q^2\alpha_n^3+q(1+q)\alpha_n^2+(2+q)\alpha_n+2-\alpha_n\alpha_k^2
\bigl[\alpha_n^3+(1-q^{-1})\alpha_n^2+(2-q)\alpha_n-2q^{-1}\bigr]\\
-\alpha_k\bigl[q^2\alpha_n^6+q(1+q)\alpha_n^5+(2+q+2q^2)\alpha_n^4+(1+q)
\alpha_n^3+2\alpha_n^2-(2+q+q^{-1})\alpha_n+(q^{-1}-1)\bigr],\endmultline$$

$y_0(n)= -q(\alpha_n-1)(\alpha_n+1)(q^2\alpha_n^2+q\alpha_n+2)$, ${ }$ 
$y_2(n)=-(q\alpha_n-1)(q\alpha_n+1)(\alpha_n^2+\alpha_n+2)$, $ {}$ 
$$
y_1(n)=q^4\alpha_n^7+q^2(1+q)(q\alpha_n^6+\alpha_n^4)+q(1+q+q^2)(2q\alpha_n^5+\
alpha_n^3)-(1+3q+3q^2+q^3)\alpha_n^2-(1+q^2)(2+\alpha_n),$$
and $\alpha_n=q^{n+1}$, $ { }$ $\beta_k=q^{k+1}.$ 

\bigskip
Then $$\alignat3
L(b(n,k))&=B(n,k)-B(n,k-1) &\qquad\text{and} &\qquad 
L(b(n,k)c(n,k))&=A(n,k)-A(n,k-1).\tag{$**$}\endalignat$$
Now, summing over k in $(**)$ shows that both sequences $a(n)$ and $b(n)$ are 
solutions of $Lu(n) = 0$.

\pagebreak
\midspace{.4in}

\bf 2.6. \rm Similar arguments and estimates as in (1.6) above lead to
$$Ln_q(2)-\frac {a_n}{b_n}=O\(b_n^{-2}\).\tag2.6.1$$
In particular, the sequence of rational numbers $\frac {a_n}{b_n}$ converges 
moderately quickly to $Ln_q(2)$.

\bf 2.7. \rm
\bf Lemma 3: \rm The sequences 
$$
v_n=a_n\prod_{t=1}^n(1+q^t)\prod_{s=[n/2]}^n(1-q^s) \qquad\text{and} \qquad 
w_n= b_n\prod_{t=1}^n(1+q^t)\prod_{s=[n/2]}^n(1-q^s)$$
are polynomials in $q$ with integer coefficients, and moreover 
$$w_n = O\(q^{19n^2/8}\).\tag2.7.1$$
  
\bf Proof: \rm Applying (1.7.1) and (1.7.2), we have estimates for the 
denominator of $v_n$: 
$$\align
ord_p\(\frac {(1-q^m)(q^2)_{n}\binom{n+m}m_q}{\binom{n+k}k_q(q)_n}\)
&\leq ord_p\(\frac {(q^m-1)(q^2)_{n}\binom km_q}{\binom{n+k}{k-m}_q(q)_n}\)\\
&\leq ord_p\(\frac{(q^2)_n}{(q)_{n}}\) + ord_p(q^m-1) + ord_p(q)_k - 
ord_p(q)_m\\
&\leq ord_p\(\frac{(q^2)_n}{(q)_{n}}\) + ord_p\prod_{s=[n/2]}^n(1-q^s) +  
ord_p(q)_k - ord_p(q)_m\\
&\leq ord_p\(\prod_{t=1}^n(1+q^t)\prod_{s=[n/2]}^n(1-q^s)\),\endalign$$
since $m \leq k \leq n$. This proves the claim on $v_n$. And (2.7.1) follows 
from (1.6.1). The rest is trivial. \bigskip
 
\bf Lemma 4: \rm $Ln_q(2)$ - $\frac {v_n}{w_n}$ = $O\(\frac 
1{w_n^{1+\delta}}\)$; where $\delta =0.26316\dots$$>$ 0.\bigskip

\bf proof: \rm Combining (1.6.1), (2.6.1) and (2.7.1), we find that
$$
Ln_q(2) - \frac{v_n}{w_n}=O\(b_n^{-2}\)
=O\(q^{-3n^2}\)
=O\(w_n^{-1-(5/19)}\).$$
\smallskip
Thus, we have proved:\smallskip
\bf Theorem 2: \rm If $|q|\neq 0, 1$ is an integer, $Ln_q(2)$ is irrational 
with irrationality measure 4.80.\smallskip
\bf Remark 2: \rm We invoke Theorem 7 ([Z], p. 596) with $\omega$ as in 2.1, 
to get the accelerated series:
$$
Ln_q(2)=\sum_{n=1}^{\infty}\frac {q^n(q)_n}{(1-q^n)(q^2)_n} 
\qquad\text{and}\qquad
Ln_q(2)=\sum_{n=1}^{\infty}\frac {(-1)^{n-1}(q)_n(1-q^{3n})}{(1-q^n)^2\binom 
{2n}n_q (q^2)_n}.$$

\pagebreak
\midspace{.3in}

\enddocument